\include{BoxedEPS}
\include{BoxedEPS}

\documentclass[12pt]{amsart}
\usepackage{amsmath,amsthm}
\input amssym.def
\input amssym

\chardef\bslash=`\\ 

\makeatletter
\def\verbatim{\interlinepenalty\@M \@verbatim
   \leftskip\@totalleftmargin\advance\leftskip2pc
   \frenchspacing\@vobeyspaces \@xverbatim}
\makeatother
\newtheorem{thm}{Theorem}[section]
\newtheorem{cor}[thm]{Corollary}
\newtheorem{lem}[thm]{Lemma}
\newtheorem{prop}[thm]{Proposition}

\theoremstyle{definition}
\newtheorem{defn}{Definition}[section]

\theoremstyle{remark}
\newtheorem{rem}{Remark}[section]
\newtheorem{exmp}{Example}[section]

\numberwithin{equation}{section}

\newcommand{\begeq}{\begin {equation}}
\newcommand{\eq}{\end{equation}}
\newcommand{\bs}{\begin {split}}
\newcommand{\es}{\end{split}}
\newcommand{\bp}{\begin {prop}}
\newcommand{\ep}{\end {prop}}
\newcommand{\bt}{\begin {thm}}
\newcommand{\et}{\end {thm}}
\newcommand{\bc}{\begin {cor}}
\newcommand{\ec}{\end {cor}}
\newcommand{\bl}{\begin {lem}}
\newcommand{\el}{\end {lem}}
\newcommand{\bpf}{\begin {proof}}
\newcommand{\epf}{\end {proof}}
\newcommand{\bi}{\begin {itemize}}
\newcommand{\ei}{\end {itemize}}
\newcommand{\ben}{\begin {enumerate}}
\newcommand{\een}{\end {enumerate}}
\newcommand{\brem}{\begin {rem}}
\newcommand{\erem}{\end {rem}}
\renewcommand{\kappa}{\wp}

\newcommand{\norm}[1]{\left\|{#1} \right\|}

\newcommand{\abs}[1]{\left|{#1} \right|}

\newcommand{\la}{\langle}
\newcommand{\ra}{\rangle}

\newcommand{\rf}[1]{\lceil {#1} \rceil}

\newcommand{\wegt}{\omega}

\newcommand {\sumzd}[1] {\sum\limits_{#1 \in \zd} }

\newcommand {\Id}{ {[0,1]^d}}

\newcommand{\AAA}{\mathbb A}
\newcommand{\BBB}{\mathbb B}

\newcommand{\HH}{{\mathcal H}}

\newcommand{\PP}{{I\kern-.3em P}}
\newcommand{\ZZ}{{Z\kern-.45em Z}}
\newcommand{\TT}{{T\kern-.45em T}}
\newcommand{\RR}{{I\kern-.3em R}}
\newcommand{\CC}{{I\kern-.6em C}}
\newcommand{\NN}{{I\kern-.3em N}}
\newcommand{\GG}{\ZZ^d}

\newcommand{\F}{\mathcal{F}}

\newcommand{\X}{\mathcal{X}}

\newcommand{\E}{\mathcal{E}}

\newcommand{\zd}{{\ZZ ^d}}

\newcommand{\Lp}{{L^p} }

\newcommand{\lp}{{\ell^p} }
\newcommand{\lqq}{{\ell^q} }

\newcommand{\linf}{{\ell^{\infty}} }
\newcommand{\lt}{{\ell^2} }
\newcommand{\betak}{{\tilde{\alpha}}}

\newcommand{\Aleph}{\aleph}

\newcommand{\OOO}{\mathcal O}
\newcommand{\ooo}{\, _{^{\mathcal O}}}
\newcommand{\ccc}{ \mathfrak c}

\newcommand{\dumvar}{\cdot}

\DeclareMathOperator{\ssup}{ess\; sup}

\DeclareMathOperator{\sprt}{supp }

\begin{document}

\bibliographystyle{plain}

\title[SLANTED MATRICES, FRAMES, AND SAMPLING] {SLANTED MATRICES, BANACH FRAMES, \\ AND SAMPLING}

\author{Akram Aldroubi, Anatoly Baskakov, and Ilya Krishtal}
\address{Department of Mathematics, Vanderbilt University, Nashville, TN 37240 \\
email: aldroubi@math.vanderbilt.edu}
\address{Department of Applied Mathematics and Mechanics, Voronezh State University, Voronezh, Russia 394693 \\
email: mmio@amm.vsu.ru}
\address{Department of Mathematical Sciences, Northern Illinois University, DeKalb, IL 60115 \\
email: krishtal@math.niu.edu}

\thanks{ The first author was supported in part by NSF grants DM-0504788.
The second author is supported in part by RFBR grant 07-01-00131.}


\date{\today }

\subjclass[2000]{Primary: 47B37, 42C15, 47N99, 94A20}

\keywords{Slanted matrices, boundedness  below, Banach frames, irregular sampling, non-uniform sampling}

\begin{abstract}
In this paper we present a rare combination of abstract results on the spectral properties
of slanted matrices and some of their very specific applications to frame theory and sampling problems.
We show that for a large class of slanted matrices boundedness below of the corresponding operator
in $\lp$ for some $p$ implies boundedness below in $\lp$ for all $p$. We use the established result
to enrich our understanding of Banach frames and obtain new results for irregular sampling problems.   
We also present a version of a non-commutative Wiener's lemma for slanted matrices.
\end{abstract}

\maketitle

\section{Introduction.}

Given a sampling set for some $p\in [1,\infty]$ can we deduce that this set is a set of sampling for all $p$?
Under which conditions a $p$-frame for some $p\in [1,\infty]$ is a Banach frame for all $p$? These and many other problems 
can be reformulated in the language of infinite slanted matrices. This paper presents such a reformulation and gives an answer to these and other interesting questions.

We begin with a simple motivational example for slanted matrices that comes from sampling theory in Paley-Wiener spaces. 
 It is well known that the Paley-Wiener space $PW_{1/2}=\{f\in L^2(\RR): \sprt \hat f \subset [-1/2,1/2]\}$ can also be described as  \begin {equation}
\label {SIS}
PW_{1/2}=\{f\in L^2(\RR): f =\sum_{k\in \ZZ}c_k\phi(\dumvar-k), \; c \in \ell^2(\ZZ)\},
\end {equation}
where $\phi(x)=\frac {\sin \pi (x -k)}{\pi (x -k)}$ and the series converges in $L^2(\RR)$ (see e.g., \cite {aldgr:siam}). Because of this equivalent description of $PW_{1/2}$, the problem of reconstructing a function $f\in PW_{1/2}$ from the sequence of its integer samples, $\{f(i)\}_{i\in \ZZ}$, is equivalent to finding the coefficients $c\in \ell^2$ such that $\{f(i)\}=Ac$ where $A=(a_{i,j})$ is the matrix with entries $a_{i,j}=\phi(i-j)$. 
It is immediate, however, that $A=I$ is the identity matrix and, therefore, 
\[f= \sum_{k\in \ZZ}f(k)\phi(\dumvar-k).\] 
If, instead, we sample a function $f \in PW_{1/2}$ on $\frac {1}{2}\ZZ$, then we obtain  the equation $\{f(\frac {i}{2})\}=Ac$. In this case, the matrix $A$ is defined by $a_{i,j}=\phi(\frac {i}{2}-j)$ and is no longer diagonal -- it has constant values on slanted lines with slopes $1/2$, for instance, $a_{2j,j}=1$.  If $\phi=\frac {\sin \pi (x -k)}{\pi (x -k)}$ in \eqref{SIS} is replaced by a function $\psi$ supported on $[-\frac {M}{2},\frac {M}{2}]$, then the matrix $A=(a_{i,j})$ is zero outside the slanted band $|j-i/2|\le M$. Again, this matrix is not banded in the classical sense. 

Another example where slanted matrices arise naturally is in the study of refinement equations and the implementation of wavelet filter algorithms (see for example \cite {BJ02} and the references therein).  
The reader should also keep in mind that the standard case of banded matrices is a particular case of slanted banded matrices. 

\subsection{Slanted Matrices: Definitions.}
\label {SMD}
\
 We prefer to give a straightforward definition
of slanted matrices in the relatively simple case that arises in applications presented
in this paper, mainly in connection with sampling theory.
For that reason, we restrict our attention to the group $\GG$, $d\in\NN$, and leave
the case of more general locally compact Abelian groups 
for future research in the spirit of \cite{Bas97, Bas97', BK05}. We believe, also, that
some of the results below may be extended to matrices indexed by discrete metric spaces.

For each $n\in\GG$ we let $X_n$ and $Y_n$ be (complex) Banach spaces and $\ell^p=\ell^p(\GG,(X_n))$ 
be the Banach space of sequences $x = (x_n)_{n\in\GG}$, $x_n\in X_n$, with the norm 
$\norm{x}_p = \left(\sum\limits_{n\in\GG}\norm{x_n}_{X_n}^p\right)^{\frac1p}$
when $p\in[1,\infty)$ and $\norm{x}_\infty = \sup\limits_{n\in\GG}\norm{x_n}_{X_n}$. 
By $\ccc_0 = \ccc_0(\GG,(X_n))$ we  denote the subspace of $\ell^\infty$ of sequences
vanishing at infinity, that is $\lim\limits_{|n|\to\infty}\|x_n\| = 0$, where $|n| = \max\limits_{1\leq k\leq d} |n_k|$,
$n = (n_1,n_2,\dots,n_d)\in\GG$.
We will use this multi-index notation throughout the paper.
Let $a_{mn}: X_n\to Y_m$ be bounded linear operators. The
symbol $\AAA$ will denote the operator matrix $(a_{mn})$, $m,n \in \GG$. In this paper, we are interested only in those
matrices that give rise to bounded linear operators that map $\ell^p$ into $\ell^p$ for all $p\in[1,\infty]$ and $\ccc_0$ into $\ccc_0$.  
We let $\norm{\AAA}_p$ be the operator norm of $\AAA$ in $\ell^p(\GG,(Y_n))$ and
$\norm{\AAA}_{\sup} = \sup\limits_{m,n\in\GG}\norm{a_{mn}}$. If $X_n$, $Y_n$, $n\in\GG$, are separable Hilbert spaces,  we denote  
by $\AAA^\star =(a_{mn}^\star)$ the matrix
defined by $a_{mn}^\star = a_{nm}^*$, where $a_{nm}^*:Y_n\to X_m$ are the
(Hilbert) adjoints of the operators $a_{nm}$.
Clearly, $(\AAA^\star)^\star = \AAA$.

\begin{rem}
Most of the applications presented in this paper are restricted to operators on $\lp(\GG,\CC)$.
The reason we introduce the spaces $\lp(\GG, (X_n))$ is twofold. First, the method we use to
prove the abstract results is not affected by this generalization. Secondly, the spaces $\lp(\GG, X)$
are better suited for applications to differential equations \cite{Bas99} and the spaces
$\lp(\GG, (X_n))$  can be used to study differential equations with unbounded operator coefficients \cite{BasKri, BasPas}.
\end{rem}

To define certain classes of operator matrices we use the following types of weight functions.

\begin{defn}
A \textit{weight} is a function $\wegt: \GG\to[1,\infty)$. A weight is \textit{submultiplicative} if
\[\wegt(m+n)\leq C \wegt(m)\wegt(n), \ \mbox { for some } C > 0.\]
A weight is a \textit{GRS}-weight if it satisfies the Gelfand--Raikov--Shilov \cite{GRS} condition
\begeq\nonumber 
\lim_{m\to \infty}\wegt(mn)^{\frac1m}=1,\ n\in\GG.
\end{equation}
A weight is \textit{balanced} if
\begeq\nonumber 
\sup_{n\in\GG}\frac {\wegt(k n)}{\wegt(n)} <\infty,\ k\in\NN.
\end{equation}
Finally, an \textit{admissible} weight is an even  submultiplicative weight.
\end{defn}

\begin{exmp}
A typical weight on $\GG$ is given by
\[\wegt(n) = e^{a|n|^b}(1 + |n|)^s, \ a,b,s\geq 0.\]
This weight is admissible when $b\in[0,1]$, is a GRS-weight when $b\in[0,1)$
and is balanced when $b=0$.
\end{exmp}

Throughout the paper we fix a slant $\alpha \neq 0$. To simplify the notation
we use  $\beta = \alpha^{-1}$ and $K = \rf{|\beta|}^d$ -- the $d$-th power of the smallest integer number bigger than or equal to $|\beta|$.
By $\chi_S$ we denote the characteristic function of a set $S$.

\begin{defn}
For $\alpha \neq 0$ and $j\in\GG$ the matrix $A_j = A_j^\alpha = (a^{(j)}_{mn})$,  $m,n\in\GG$, defined by 
\[a^{(j)}_{mn} = a_{mn}\prod_{k=1}^d\chi_{[j_k,j_k+1)}(\alpha m_k- n_k)\] is called the
{\it $j$-th $\alpha$-slant} of $\AAA$.
\end{defn}

Observe that for every $m\in\GG$ there is at most one $n\in\GG$ such that $a^{(j)}_{mn}\neq 0$
and at most $K$ different numbers $\ell\in\GG$ such that $a^{(j)}_{\ell m}\neq 0$.
Hence, we have $\|A_j\|_p \leq K\|A_j\|_{\sup}$
for any $p\in[1,\infty]$. This allows
us to define different classes of matrices with decaying $\alpha$-slants
independently of $p\in[1,\infty]$. 

\begin{defn}\label{slantm} 
We consider the following several classes of matrices.
\begin{itemize}
\item For some fixed $M\in\NN$, 
$\F_\alpha^M$ will denote the class of matrices $\AAA$ that satisfy $\AAA = \sum_{|j|\leq M-1} A_j$. 
Observe that for $\AAA\in\F_\alpha^M$ we have $a_{mn} = 0$ as soon as $|n-\alpha m| > M-1$.
The class $\F_\alpha = \bigcup_{M\in\NN} \F_\alpha^M$ 
consists of operators with {\it finitely many $\alpha$-slants}.

\item The class $\Sigma^\wegt_\alpha$ of matrices with {\it $\wegt$-summable} $\alpha$-slants 
consists of matrices $\AAA$  
such that $\norm{\AAA}_{\Sigma^\wegt_\alpha} = K\sum_{j\in\GG}\norm{A_j}_{\sup}\wegt(j) <\infty$, where $\wegt$ is a  weight. We have $\Sigma^\wegt_\alpha \subset\Sigma^1_\alpha = \Sigma_\alpha$  -- the class of matrices with (unweighted) summable $\alpha$-slants. 

 \item The class $\E_\alpha$ of matrices with {\it exponential decay} of $\alpha$-slants is defined as a subclass
of matrices $\AAA$ from $\Sigma_\alpha$ such that for some $C\in\RR$ and $\tau\in(0,1)$
we have $\norm{A_j}_{\Sigma_\alpha} \leq C\tau^{|j|}$.
\end{itemize}
\end{defn}

  For $\AAA\in\Sigma^\wegt_\alpha$, we  denote by $\AAA_M\in \F_\alpha^M$, $M\in\NN$, the {\it truncation} of $\AAA$, {\it i.e.}, the matrix defined by
$a_{mn}^M = a_{mn}$ when $|n-\alpha m| \leq M-1$ and $a_{mn}^M = 0$ otherwise. Equivalently, $\AAA_M = \sum\limits_{|j|\leq M-1} A_j$, where $A_j$, $j\in\GG$, is the $j$-th $\alpha$-slant of $\AAA$. By definition of $\Sigma^\wegt_\alpha$, 
the operators $\AAA_M$ converge to $\AAA$ in the norm  $\norm{\cdot}_{\Sigma^\wegt_\alpha}$.

\begin{rem}\label{R12}
Notice that when $\alpha=d=1$ we get the usual matrix diagonals as a special case of $\alpha$-slants studied
in this paper.  
We introduce $\alpha$-slants, in part, to avoid certain reindexing problems that occur in applications. These
kind of problems are treated differently in \cite{BCHL, Qiyu}. In particular, it is not hard to see that in many cases a slanted matrix can be converted into a conventional banded block
matrix.
\end{rem}

\begin{rem}
\label {R22}
The observation preceding Definition \ref{slantm} implies
that matrices in $\F_\alpha$ define bounded operators in $\ccc_0$,
and any $\ell^p$, $p\in[1,\infty]$. Since for any $\AAA\in\Sigma^\wegt_\alpha$
we have $\AAA = \sum_{j\in\GG}{A_j}$, where the series converges in the norm of $\Sigma^\wegt_\alpha$,
the matrices in 
$\Sigma^\wegt_\alpha$ and $\E_\alpha$ also define bounded operators in $\ccc_0$,
and any $\ell^p$, $p\in[1,\infty]$. Moreover, it is not hard to see that
$\Sigma^\wegt_\alpha$ is a Banach space with respect to the norm given by $\norm{\cdot}_{\Sigma^\wegt_\alpha}$ and
$\norm{A}_p \leq \norm{A}_{\Sigma_\alpha}\leq \norm{A}_{\Sigma^\wegt_\alpha}$ for every $A\in \Sigma^\wegt_\alpha$ and $p\in[1,\infty]$.
We also point out that all results in Section \ref{pure}
that are stated for $\ell^\infty$ also hold for $\ccc_0$, although we do not mention it later.
\end{rem}

\begin {rem} \label {OpToMat} There is an obvious one to one and onto correspondence between the matrices in $\Sigma^\wegt_\alpha$ and a class of operators on $\ell^p$. In particular, given an operator $\BBB: \ell^p(\GG,(X_n))\to \ell^p(\GG,(Y_n))$, we can define its matrix $b_{mn}:X_n\to Y_m$ by restricting the operator $P^Y_m\BBB P^X_n$, where $P^X_n$   are given by 
\[
P^X_n(\cdots,x_{n-1},x_n,x_{n+1},\cdots)=(\cdots,0,x_n,0,\cdots),\ n\in\GG,
\]
and $P^Y_m, \; m\in\GG$, are defined in a similar way. Below, we do not distinguish between an operator and its matrix when no confusion may arise.

Moreover, one can define a matrix of an operator on any Banach space given a \textit{resolution of the identity}, which is a family of projections
with the same properties as $P^X_n$, $n\in\GG$. We refer to \cite{Bas97, Bas97'} for details.
\end {rem}

	The rest of this paper is organized as follows. We study some basic properties of  slanted matrices in \ref {SMBP}. In Section \ref {MPR} we state and prove one of our main theorems. Specifically, slanted matrices with some decay, viewed as operators on $\ell^p$ spaces,  are either universally bounded below for all $p\in [1,\infty]$, or do not have this property for any $p \in [1,\infty]$. In Section \ref {WTLSC} we use this theorem to obtain a version of Wiener's Tauberian Lemma and a result on subspace complementation in Banach spaces. Section \ref {A} is devoted to some applications of the results of Section \ref {AR}. Specifically, in Section \ref {BF}, the reconstruction formulae for Hilbert frames are extended to Banach frames under certain localization conditions related to slanted matrices. Gabor systems having this localization property are then presented as an example. Section   \ref {SRP} exhibits an application of slanted matrices to sampling theory. 
\section{Abstract Results}\label{pure}
\label {AR}

\subsection{Slanted Matrices: Basic Properties.}
\label {SMBP}

Here we present some basic properties of slanted matrices that are useful for the remainder of the paper. 
\bl
\label{bp}
For some $p\in[1,\infty]$ we consider two operators  $\AAA:\ell^p(\GG,(Y_n))\to\ell^p(\GG,(Z_n))$ and
$\mathbb{B}:\ell^p(\GG,(X_n))\to\ell^p(\GG,(Y_n))$  and let $\wegt$ be a submultiplicative balanced weight. 
\begin{itemize}
\item If $\AAA\in\F_\alpha\,(\Sigma^\wegt_\alpha,\,or\, \E_\alpha)$ and $\mathbb{B} \in\F_\betak\,(\Sigma^\wegt_\betak,\,or\, \E_\betak)$
then we have $\AAA\mathbb{B}\in\F_{\alpha\betak}\,(\Sigma^\wegt_{\alpha\betak},\,or\, \E_{\alpha\betak})$.
\end{itemize}
If, moreover, $Y_n$, $Z_n$, $n\in\GG$, are Hilbert spaces, then we have $\AAA^\star:\ell^p(\GG,(Z_n))\to\ell^p(\GG,(Y_n))$ and
\begin{itemize}
\item $\AAA$ is invertible if and only if $\AAA^\star$ is invertible;
\item If $\AAA\in\F_\alpha\,(\Sigma^\wegt_\alpha,\,or\, \E_\alpha)$ then $\AAA^\star\in
\F_{\alpha^{-1}}\,(\Sigma^\wegt_{\alpha^{-1}},\,or\, \E_{\alpha^{-1}})$;
\end{itemize}
\el
\bpf
The last two properties are easily verified by direct computation. For the first one, let
 $\mathbb {D}=(d_{m,n})=\AAA\mathbb{B}=(a_{m,n})(b_{m,n})$, and let $\rf{a} = (\rf{a_1},\dots, \rf{a_d})\in\GG$, where $a\in\RR^d$ and $\rf{a_k}$ is, as before, the smallest integer greater than or equal to $a_k$, $k = 1,\dots, d$. We have that  
\begeq\nonumber
\begin{split}
\norm {d_{m, \rf {\alpha\betak m}+j}}&\le \sum_{k\in\GG} \norm {a_{m,k}} \norm {b_{k,\rf {\alpha\betak m}+j}} \\
&=\sum_{k\in\GG} \norm {a_{m,\rf{\alpha m}+k-\rf{\alpha m}}} \norm {b_{k,\rf {\betak k} +\rf {\alpha\betak m}+j - \rf {\betak k}}}\\
&\le \sum_{k\in\GG} r\big({k-\rf{\alpha m}}\big)s\big({\rf {\alpha\betak m}+j - \rf {\betak k}}\big)\\
&= \sum_{k\in\GG} r\big(k\big)s\big({\rf {\alpha\betak m}+j - \rf {\betak k+\betak \rf{\alpha m} }}\big),
\end {split}
\end{equation}
where $r(j)=\norm {A_j}_{\sup}$ and $s(j)=\norm {B_j}_{\sup}$. For $a, b\in\RR$ we have
$$\rf {a}+\rf{b}-1\le \rf{a+b}\le \rf {a}+\rf{b};$$ 
$$
\rf{|a|b}\le \rf {|a|\rf{b}}\le\rf{|a|b}+\rf{|a|}.
$$
Hence, 
$$
\norm {d_{m, \rf {\alpha\betak m}+j}}\le  \sum_{k\in\GG} r\big(k\big)s\big({j - \rf {\betak k }+l}\big),
$$
where $l=l(\alpha,\betak,m,k)\in\GG$ is  such that $|l|\le \rf {|\betak|}+1$. 

If $\AAA\in\F_\alpha$ and $\mathbb B\in \F_\betak$, the last inequality immediately implies 
$\mathbb D = \AAA\mathbb B\in \F_{\alpha\betak}$.

If $\AAA\in\Sigma^\wegt_\alpha$ and $\mathbb B\in \Sigma^\wegt_\betak$, we use the fact that the weight
$\wegt$ is submultiplicative and balanced to obtain
\[\sum_{j\in\GG}\sup_{m\in\GG}\norm {d_{m, \rf {\alpha\betak m}+j}}\wegt(j)\le  \sum_{j,k\in\GG} r\big(k\big)
s\big({j - \rf {\betak k }+l}\big)\wegt(j)\le\]
\[Const\cdot\sum_{j,k\in\GG} r\big(k\big)\wegt\big(k\big)s\big({j - \rf {\betak k }}\big)\wegt(j - \rf{\betak k})\frac{\wegt\big(\rf{\betak k}\big)}{\wegt\big(k\big)}\]
\[\le Const\cdot \norm{\AAA}_{\Sigma^{\wegt}_{\alpha}}\norm{\mathbb B}_{\Sigma^{\wegt}_{\betak}}.\]

The case $\AAA\in\E_\alpha$ and $\mathbb B\in \E_\betak$ can be treated in a similar way. Since we will not use this result in the paper,
we omit the proof. \epf

\subsection{Main Result.}
\label {MPR}
\
The property of (left, right) invertibility of operator matrices in certain operator
algebras has been studied extensively by many authors. The main focus in this paper, however, is on a
weaker property of boundedness below (or uniform injectivity). As we show in Section \ref{A} matrices with this property play a crucial role in certain applications.

\begin{defn}
We say that the matrix $\AAA$ is {\it bounded below in} $\ell^p$ or, shorter, {\it $p$-bb}, if
\begeq\label{coerp}
\norm{\AAA x}_p \geq \kappa_p \norm{x}_p, \mbox{ for some } \kappa_p > 0 \mbox{ and all } x\in\lp.
\end{equation}
\end{defn}

Before we state our main result, we note an important spectral property of slanted matrices given by the
following lemma due to Pfander \cite {P06} (see also \cite{PW02}). We include the proof for completeness and since the matrices considered here are more general.

\bl
 \label {PL}
Assume that $X_n = Y_n$, $n\in\GG$, and that all these spaces are finite dimensional. 
If $\AAA \in \Sigma_\alpha$, for some $\alpha>1$, then $0$ is an approximate eigen-value of $\AAA:\ell^p\to\ell^p$, $p\in[1,\infty]$.
Equivalently, for any $\epsilon>0$ there exists $x\in\ell^p$ such that $\norm{x}_p=1$ and $\norm{\AAA x}_p\leq\epsilon$. 
\el 

\bpf
Let $\AAA \in \Sigma_\alpha$. For  $\epsilon>0$ choose $M$ so large that 
$\norm {\AAA-\AAA_M} _{\Sigma_\alpha}\le \epsilon$. Since $\alpha>1$, there exists $N_0$ such that $N=\rf {\alpha N_0}\ge N_0+1$. Let $\AAA^N_M$ be a matrix with an $(i,j)$-entry coinciding with that of the truncation matrix $\AAA_M$ if $|i|\le M+N$, $|j|\le M+N$, and equal to $0$ otherwise. 
We have $\AAA_Mx^N_M=\AAA^N_Mx^N_M$ for every $x^N_M\in \ell^p$ such that $x^N_M(i)=0$ for $|i|>M+N$. By assumption, the subspace $\X_M^N$ of such vectors is
finite dimensional and, by construction, it is invariant with respect to $\AAA^N_M$. Observe that we chose $N$ so large that the restriction of $A_M^N$ to $\X_M^N$ cannot be invertible because its matrix has a zero ``row''. 
Hence, for $\AAA^N_M$, we can find a vector $x^N_M\in\X^N_M$ such that $\norm {x^N_M}=1$ and $\AAA_Mx^N_M=\AAA^N_Mx^N_M=0$. Thus, for any given $\epsilon>0$, we can find  $x^N_M\in\X$ such that $\norm {x^N_M}=1$, and $\norm {\AAA x^N_M}_p=\norm {\AAA x^N_M-\AAA_M x^N_M}_p\le \epsilon$.
\epf

The following theorem presents our central theoretical result.

\bt\label{mainthm}
Let $s>(d+1)^2$ and $\wegt=(1+|j|)^s$. Then
$\AAA\in\Sigma^\wegt_\alpha$ is $p$-bb for some $p\in[1,\infty]$ if and only if $\AAA$ is $q$-bb for all $q\in[1,\infty]$. 
\et

\begin{rem}\label{R123}
Observe that if  $X_n = Y_n$, $n\in\GG$, 
Lemma \ref{bp} allows us to consider only the case $\alpha > 0$. Indeed, if $\AAA\in\Sigma^\wegt_\alpha$ then
$\AAA^2\in\Sigma^\wegt_{\alpha^2}$ and it is immediate that $\AAA$ is $p$-bb if and only if $\AAA^2$ is $p$-bb.
Another way to see that we can disregard the case $\alpha < 0$ (even when $X_n \neq Y_n$) follows from the fact that the lower bound $\kappa_p$
does not change when we permute the ``rows'' of the matrix $\AAA$. Indeed, if $\AAA \in \Sigma^\wegt_\alpha$, let $\BBB = (b_{i,j})$ be defined
by $b_{i,j} = a_{-i,j}$. Then $\BBB\in\Sigma^\wegt_{-\alpha}$ and $\AAA$ is $p$-bb if and only if $\BBB$ is $p$-bb.

Observe, also, that Lemma \ref{PL} implies that often Theorem \ref{mainthm} is vacuous for $\alpha > 1$. 
Remark \ref{R12}, on the other hand, indicates that sometimes the theorem can be reduced to the case $\alpha = 1$. However, we find such a reduction 
misleading. Firstly, it does not significantly simplify our proof and, secondly, it can make computing explicit estimates
in applications more complicated.
\end{rem}

The proof of the theorem is preceded by several technical lemmas and observations below. We begin
with a lemma that provides some insight into the intuition behind the proof. We should also mention that
our approach is somewhat similar to Sj\"ostrand's proof of a non-commutative Wiener's lemma \cite{Sjo}. 
We will discuss Wiener-type lemmas in more detail in the next section.

Let  $w^N:\RR^d\to\RR$, $N>1$, be a family of window functions such that $0\le w^N\le1$, $w^N(k)=0$ for all $|k|\ge N$, and $w^N(0)=1$. 
By $w_n^N$ we will denote the translates of $w^N$, {\it i.e.}, $w^N_n(t) = w^N(t-n)$, and 
$W^N_n: \ell^p(\GG,X)\to\ell^p(\GG,X)$ will be the multiplication operator 
\[W_n^Nx(k) = w_n^N(k)x(k),\ \ x\in\ell^p,\ n\in\RR^d.\]
Let $x\in \ell^p(\ZZ^d,X)$, $p \in [1,\infty]$, and define 
\begeq \nonumber 
\bs
|||x|||_p^p&:=\sum_{n\in\ZZ^d}\sum\limits_{j\in \ZZ^d}{\|W_n^Nx(j)\|^p_p}=\sum_{n\in\ZZ^d}\sum\limits_{|j-n|< N}{\|W_n^Nx(j)\|^p_p}, \  p\in[1,\infty),\\
|||x|||_\infty&:=\sup\limits_n\|W^N_nx\|_\infty.
\end {split}
\end {equation}

\bl \label {normeq}For any $p\in [1,\infty]$, the norms $\|\dumvar\|_p$ and $|||\dumvar |||_p$ are equivalent norms on $\ell^p$, and we have 
\[
\|x\|_p\le|||x|||_p\le (2N)^{d/p}\|x\|_p, \quad p \in [1,\infty),
\]
and
\[
\|x\|_\infty=|||x|||_\infty.
\]
\el
\bpf
For $p=\infty$ the equality is obvious. For $p\in [1,\infty)$, the left inequality follows from the fact $\|x(n)\|^p\le \sum\limits_{|j-n|\le N}{\|W_n^Nx(j)\|^p_p}$, and by summing over $n$. For the right inequality we simply note that
\[
\sum_{n\in\ZZ^d}\sum\limits_{|j-n|< N}{\|w_n^N(j)x(j)\|^p_p}\le \sum\limits_{n\in\ZZ^d}\sum_{|j|< N}{\|x(j+n)\|^p_p}
\le (2N)^d\|x\|_p^p.
\]
\epf

The above equivalence of norms will supply us with the crucial inequality in the proof of the theorem. The opposite
inequality is due to the following observation.

\begin{rem}\label{find}
We shall make use of the following obvious relation between the norms in finite dimensional spaces.
For every $x$ in a $d$-dimensional Euclidean space we have
\begeq\label{pl2}
\norm{x}_p \geq \norm{x}_\infty \geq d^{-\frac1p}\norm{x}_p \mbox{ for any } p \in [1,\infty).
\end{equation}
\end{rem}

At this point we choose our window functions to be the family of Ces\`aro means
$\psi^N:\RR^d\to \RR$, $N>1$, defined by
\begeq\nonumber 
\psi^N (k) = \left\{
\begin{array}{rl}
(1-\frac{|k|}N), & |k| < N;\\
0,&\mbox{ otherwise.}
\end{array}
\right.
\end{equation}
Observe that their translates $\psi^N_n (k) = \psi^N(k-n)$, $n\in\RR^d$,
satisfy 
\begeq\label{scale}
\psi^{\alpha N}_{\alpha n} (k) = \psi^N_n (\alpha^{-1}k)
\end{equation}
 for any $\alpha > 0$.
Again, by $\Psi_n^N: \ell^p\to \ell^p$, $N>1$, we will denote the operator of multiplication \[\Psi_n^Nx(k) = \psi_n^N(k)x(k),\ \ x\in\ell^p,\ n\in\RR^d.\]

The following lemma presents yet another estimate crucial for our proof. We remind that to simplify the
notation we let $\beta=\alpha^{-1}$.

\bl\label{Cesal}
The following estimate holds for any $q\in[1,\infty]$, any $\AAA\in \Sigma_\alpha = \Sigma^1_\alpha$, and
all of its truncations $\AAA_M\in \F_\alpha^M$, $M\in\NN$.
\begeq\label{Cesest}
\norm{\AAA_M\Psi_n^N - \Psi_{\beta n}^{\beta N}\AAA_M}_q \leq \frac {(2M)^{d+1}}{2N} \norm{\AAA}_{\sup}
=:\Aleph/2. 
\end{equation}
\el

\bpf
Define $J_k=\{i \in \ZZ^d: |i-\alpha k|\le M-1\}$. Using \eqref{scale}, we have
\[\left\vert \psi_n^N(i)-\psi_{ \beta n}^{\beta N}(k) \right\vert
\leq \frac {M-1}{N}, \quad \hbox { for }  |i-\alpha k|\leq M-1.
\]
Observe that for any $y\in\ell^q$ we have
$$(\AAA_M\Psi_n^N y)(k)\! =\! \sum_{i\in J_k} a_{ki}\psi_n^N(i)y(i),\ (\Psi_{\beta n}^{\beta N} \AAA_M y)(k) = \psi_{\beta n}^{\beta N}(k)\sum_{i\in J_k} a_{ki}y(i).$$
Now the following easy computation 
\begin{equation}\nonumber
\bs
\norm{\left(\AAA_M\Psi_n^N - \Psi_{\beta n}^{\beta N}\AAA_M\right)y}_q = \left(\sum_{k\in\GG}\norm{\sum_{i\in J_k} 
a_{ki}(\psi_n^N(i)-\psi_{\beta n}^{\beta N}(k))y(i)}^q\right)^{\frac1q} \\
\leq\frac{M}{N}\norm{\AAA}_{\sup}\left(\sum_{k\in\GG}\left(\sum_{i\in J_k} 
\norm{y(i)}\right)^q\right)^{\frac1q}\leq
\frac {(2M)^{d+1}}{2N}\norm{\AAA}_{\sup}\norm{y}_q 
\end{split}
\end{equation}
shows that (\ref{Cesest}) is true for $q\in[1,\infty)$. An obvious modification yields it in the case $q=\infty$.
\epf

Observe that for $\AAA_M\in \F_\alpha^M$ the commutator studied in the above lemma satisfies
\begeq\label{crop}
(\Psi_{\beta n}^{\beta N}\AAA_M-\AAA_M \Psi_{n}^{N}) x = 
(\Psi_{\beta n}^{\beta N}\AAA_M-\AAA_M \Psi_{n}^{N}) P_n^{N+M} x, 
\end{equation}
where $\beta=\alpha^{-1}$, $P_n^L x(k) = x(k)$ if $|k-n|\leq L$, and $P_n^L x(k) = 0$ otherwise, $L>1$. 
Also observe that for any $p\in [1,\infty]$ and any $L>1$, we have that 
\begeq\label{crop2}
\norm {P_n^{L} x}_p \le 2 \norm {\Psi_{n}^{2L} x}_p.
\end{equation}

Combining the above facts we obtain the following estimate.

\bl
Let $\AAA\in \Sigma_\alpha$ be $p$-bb for some $p\in[1,\infty]$. As usually, let $\AAA_M\in \F_\alpha^M$ be the truncations of $\AAA$ and $\beta=\alpha^{-1}$.
Then for all  $n\in\GG$, $N>1$, and such $M\in\NN$ that $\gamma_p = \kappa_p - \norm{\AAA-\AAA_M}_p > 0$, we have
\begeq\label{m1}
\norm{\Psi^N_nx}_p \leq \gamma_p^{-1}\left(\norm{\Psi_{\beta n}^{\beta N }\AAA_M x}_p+\Aleph\norm{\Psi_n ^{2(N+M)}x}_p\right).
\end{equation}
\el
\bpf 
Observe that
\[\norm{\Psi_n^N x}_p \leq \kappa_p^{-1}\left(
 \norm{(\AAA_M \Psi_{n}^{N}) x}_p + \norm{\AAA-\AAA_M}_p\norm{\Psi_n^N x}_p\right).\]
Hence, using \eqref{Cesest}, \eqref{crop}, and \eqref {crop2}, we get
\begeq\nonumber
\bs
&\left(1 - \kappa_p^{-1}\norm{\AAA-\AAA_M}_p\right)\norm{\Psi_n^N x}_p   \leq \kappa_p^{-1}\norm{\AAA_M\Psi_n^N x}_p \\
&\leq \kappa_p^{-1}\left(\norm{\Psi_{\beta n}^{\beta N}\AAA_M x}_p +
 \norm{(\Psi_{\beta n}^{\beta N}\AAA_M-\AAA_M \Psi_{n}^{N}) x}_p \right)\\
&\leq \kappa_p^{-1}\left(\norm{\Psi_{\beta n}^{\beta N}\AAA_M x}_p+
 \norm{(\Psi_{\beta n}^{\beta N}\AAA_M-\AAA_M \Psi_{n}^{N}) P^{N+M}_n x}_p\right)\\
&\leq \kappa_p^{-1}\left(\norm{\Psi_{\beta n}^{\beta N}\AAA_M x}_p+ \frac{\Aleph}{2}\norm{ P^{N+M}_n x}_p\right)\\
&\leq \kappa_p^{-1}\left(\norm{\Psi_{\beta n}^{\beta N}\AAA_M x}_p+ \Aleph\norm{ \Psi^{2(N+M)}_n x}_p\right),
\end{split}
\end{equation}
which yields the desired inequality.
\epf
By iterating \eqref {m1} $j-1$ times we get
\bl\label{lemmaZ}
Let $\AAA\in \Sigma_\alpha$ be $p$-bb for some $p\in[1,\infty]$. Let $\AAA_M\in \F_\alpha^M$ be the truncations of $\AAA$ and $\beta=\alpha^{-1}$.
Then for all  $n\in\GG$, $N>1$, and such $M\in\NN$ that $\gamma_p = \kappa_p - \norm{\AAA-\AAA_M}_p > 0$, we have
\begeq\label{mj}
\norm{\Psi^N_nx}_p \leq \gamma_p^{-1}\frac {1-(\Aleph \gamma_p^{-1})^j}{1-(\Aleph \gamma_p^{-1})} 
\norm{\Psi_{\beta n}^{\beta Z_j }\AAA_M x}_p+(\Aleph \gamma_p^{-1})^j\norm{\Psi_n ^{Z_{j+1}}x}_p, 
\end{equation}
where $Z_j=2^{j-1}N+(2^j-2)M$, for $j\ge1$. 
\el

To simplify the use of \eqref{mj} we let
\begeq
\label{ajp}
a_{j,p}:=\gamma_p^{-1}\frac {1-(\Aleph \gamma_p^{-1})^j}{1-(\Aleph \gamma_p^{-1})} = 
\frac {1-\left({\frac {(2M)^{d+1}\norm{\AAA}_{\sup}}{\left(\kappa_p - \norm{\AAA-\AAA_M}_p\right) N} } \right)^j}
{\kappa_p - \norm{\AAA-\AAA_M}_p -{\frac {(2M)^{d+1}}{N} \norm{\AAA}_{\sup}}} 
\end{equation}
and 

\begeq
\label {bjp}
b_{j,p}:=(\Aleph \gamma_p^{-1})^j = 
{\frac {((2M)^{d+1} \norm{\AAA}_{\sup})^j}{\left(\kappa_p - \norm{\AAA-\AAA_M}_p\right)^j N^j}} .
\end {equation}
 
Now we are ready to complete the proof
of the main result.

\bpf (Theorem \ref{mainthm}). The remainder of the proof will be presented in two major steps. In the first step, we will
assume that $\AAA\in\Sigma^\wegt_\alpha$ is $\infty$-bb and show that this implies that $\AAA$ is $p$-bb for any $p\in[1,\infty)$.
In the second step we will do the ``opposite'', that is, assume  
that $\AAA\in\Sigma^\wegt_\alpha$ is $p$-bb for some $p\in[1,\infty)$ and show that this implies that $\AAA$ is $\infty$-bb.
This would obviously be enough to complete the proof. 

\textit{Step 1}. Assume that $\AAA$ is $\infty$-bb. Using H\"older's inequality and \eqref {mj}, we get for large values of $M\in\NN$
\begeq\nonumber
\norm{\Psi_n^N x}^p_\infty  \leq 2^{p-1}a^p_{j,\infty} \norm{\Psi_{\beta n}^{\beta Z_j }\AAA_M x}^p_\infty+2^{p-1}b^p_{j,\infty}\norm{\Psi_n ^{Z_{j+1}}x}^p_\infty.
\end{equation}
Using \eqref {pl2}, we get
\begeq\nonumber
(2N)^{-d}\norm{\Psi_n^N x}^p_p  \leq 2^{p-1}a^p_{j,\infty} \norm{\Psi_{\beta n}^{\beta Z_j }\AAA_M x}^p_p+2^{p-1}b^p_{j,\infty}\norm{\Psi_n ^{Z_{j+1}}x}^p_p.
\end{equation}
Summing over $n$ and using Lemma \ref {normeq}, we get
\begin {equation}\label{longo}
\bs
\norm{x}^p_p&\le (2N)^d2^{p-1}a_{j,\infty}^p(2 Z_j)^d\norm {\AAA_M x}_p^p+(2N)^d2^{p-1}b_{j,\infty}^p(2 Z_{j+1})^d\norm {x}^p_p\\
&\le N^d2^{2d+p-1}a_{j,\infty}^pZ_j^d\left(\norm {\AAA x}_p^p + 
\norm {\AAA - \AAA_M}_p^p \norm{x}_p^p\right) \\ &+N^d2^{2d+p-1}b_{j,\infty}^pZ_{j+1}^d\norm {x}^p_p.
\end{split}
\end {equation}

At this point we use the assumption $\AAA\in \Sigma_\alpha^{(1+|j|)^s}$ to get
\[\norm{\AAA - \AAA_M}_p \le \norm{\sum_{|j|\ge M} A_j}_p \le K\sum_{|j|\ge M} \norm{A_j}_{\sup}(1+|j|)^s (1+|j|)^{-s}\]
\[\le \norm{\AAA}_{\Sigma_\alpha^{(1+|j|)^s}}\sup_{|j|\ge M}(1+|j|)^{-s}
\le \norm{\AAA}_{\Sigma_\alpha^{(1+|j|)^s}}M^{-s}.\]

Plugging the above estimate into \eqref{longo} we obtain 
\begin {equation}\label{longo1}
\bs
\norm{x}^p_p&\le 2^{2d+p-1}N^da_{j,\infty}^pZ_j^d\norm {\AAA x}_p^p  
\\ &+2^{2d+p-1}N^d\left(a_{j,\infty}^pZ_j^d\norm{\AAA}_{\Sigma_\alpha^{(1+|j|)^s}}^p M^{-sp} + b_{j,\infty}^pZ_{j+1}^d\right)\norm {x}^p_p \\
& = 2^{2d+p-1}a_{j,\infty}^pN^dZ_j^d\norm {\AAA x}_p^p + \widetilde{\Aleph}\norm{x}^p_p.
\end{split}
\end {equation}

Hence, to complete Step 1 it suffices to show that one can choose $j, M\in \NN$ and $N > 1$ so that $\widetilde{\Aleph}< 1$. 

We put $N=M^{\delta (d+1)}$ for some $\delta > 1$. From \eqref{Cesest}, \eqref{ajp}, \eqref{bjp}, and the definition of $Z_j$
in Lemma \ref{lemmaZ} we get $\Aleph = \OOO(M^{(1-\delta)(d+1)})$, $b_{j,\infty} = \OOO(M^{(1-\delta)(d+1)j})$, 
$a_{j,\infty} = \OOO(1)$, and $Z_{j}= \OOO(M^{\delta (d+1)})$ as $M\to\infty$. Hence,
\[\widetilde{\Aleph} \le C_1 M^{\delta (d+1)^2-sp} + C_2 M^{(1-\delta)(d+1)jp + \delta (d+1)^2 }, \]
where the constants $C_1$ and $C_2$ depend on $\AAA$, $s$, $j$, and $p$ but do not depend on $M$.
Since $s > (d+1)^2$, we can choose $\delta\in(1,\frac{sp}{(d+1)^2})$ and $j>\frac{\delta (d+1)}{p(\delta - 1)}$.
Then, clearly, $\widetilde{\Aleph} = \ooo(1)$ as $M\to\infty$.

\textit{Step 2}. Now assume that $\AAA$ is $p$-bb, for some $p\in [1,\infty)$. Using  \eqref {pl2} and \eqref {mj},
we get
\begin {equation}\nonumber
\norm{\Psi_n^N x}_\infty \leq a_{j,p} (2 Z_j )^{d/p}\norm{\Psi_{\beta n}^{\beta Z_j }\AAA_M x}_\infty+b_{j,p}(2Z_{j+1})^{d/p}\norm{\Psi_n ^{Z_{j+1}}x}_\infty.
\end{equation}
As in Step 1, we have 
$\norm{\AAA - \AAA_M}_p \le \norm{\AAA}_{\Sigma_\alpha^{(1+|j|)^s}}M^{-s}$. Using this estimate
and Lemma \ref {normeq}, we obtain 
\begin {equation}\nonumber
\bs
\norm{x}_\infty &\leq a_{j,p} (2 Z_j )^{d/p}\norm{\AAA x}_\infty \\
&+ 2^{d/p}\left( a_{j,p} 
Z_j^{d/p}\norm{\AAA}_{\Sigma_\alpha^{(1+|j|)^s}}M^{-s}+b_{j,p}Z^{d/p}_{j+1}\right)\norm{x}_\infty.
\end{split}
\end{equation}
Again, as in the previous step, if we choose $\delta\in(1,\frac{sp}{(d+1)^2})$, $N = M^{\delta(d+1)}$, and $j>\frac{\delta (d+1)}{p(\delta - 1)}$, we get
\[a_{j,p} 
Z_j^{d/p}\norm{\AAA}_{\Sigma_\alpha^{(1+|j|)^s}}M^{-s}+b_{j,p}Z^{d/p}_{j+1}
= \ooo(1)\]
as $M\to\infty$ and the proof is complete.
\epf

Careful examination of \eqref{longo1} yields the following result.

\bc
Let $s>(d+1)^2$, $\wegt=(1+|j|)^s$,  and
$\AAA\in\Sigma^\wegt_\alpha$ be $p$-bb for some $p\in[1,\infty]$. Then there exists $\kappa > 0$ such that for all $q\in[1,\infty]$
\[\norm{\AAA x}_q \geq \kappa \norm{x}_q, \mbox{ for all } x\in\ell^q.\]
\ec

As we have seen in the proof above, the group structure of the index set $\GG$ has not been used. Thus, it is natural to conjecture
that a similar result holds for matrices indexed by much more general (discrete) metric spaces. In this paper, however, we do not pursue this extension.
Instead, we prove the result for a class of matrices that define  operators of \emph{bounded flow}.

\begin{defn}
A matrix $\AAA$ is said to have \emph{bounded dispersion} if there exists $M\in\NN$ such that for every $m\in\GG$ there exists $n_m\in\GG$ 
for which $a_{mn} = 0$ as soon as $|n-n_m| > M$. A matrix $\AAA$ is said to have \emph{bounded accumulation} if $\AAA^\star$ has bounded dispersion.
Finally, $\AAA$ is a \emph{bounded flow} matrix if it has both bounded dispersion and bounded accumulation.
\end{defn}

\bc
Assume that $\AAA$ has bounded flow and is $p$-bb for some $p\in[1,\infty]$. Then $\AAA$ is $q$-bb for all $q\in[1,\infty]$.
\ec

\bpf
In lieu of the proof it is enough to make the following two observations. First, if a matrix is bounded below then any matrix obtained from the original one 
by permuting its rows (or columns) is also bounded below with the same bound. Second, if a matrix is bounded below then any matrix obtained from the original one 
by inserting any number of rows consisting entirely of $0$ entries is also bounded below with the same bound. Using these observations we can use row permutations and  insertions of zero rows to obtain a slanted matrix in $\mathcal F^M_\alpha$ for some $\alpha \in \RR$, $|\alpha|>0$. 
 \epf

\brem \label{lb} The proof of Theorem \ref {mainthm} indicates how an explicit bound $\kappa_q$ and a universal bound $\kappa$ can be obtained in terms of $\kappa_p$. We did not compute these bounds because such calculations may be easier and yield better results in specific examples.
\erem

\subsection{Wiener-type Lemma and Subspace Complementation.}\ 
\label {WTLSC}
The classical Wiener's Lemma \cite{wiener} states that if a periodic function $f$ has an absolutely convergent Fourier series and never vanishes
then the function $1/f$ also has an absolutely convergent Fourier series. This result has many extensions (see \cite{Bal, Bas90, Bas97, Bas97', GKW, GroL, Jaf, Ku90, Loom, Sjo, Qiyu, Shu} and references therein), 
some of which have been used
recently in the study of localized frames \cite{BCHL, Gro}. Most of the papers just cited show how Wiener's result can be viewed as a statement
about the off-diagonal decay of matrices and their inverses. 
Using Lemma \ref {bp} and  \cite[Theorem 2]{Bas97'} we obtain the following result about invertible slanted matrices.
\bt\label{Inv}
Let $X_n$, $Y_n$, $n\in\GG$, be Hilbert spaces and $\wegt$ be an admissible balanced GRS-weight. 
If  $\AAA\in\Sigma^\wegt_{\alpha}$ is invertible for some $p\in[1,\infty]$, then $\AAA$ is invertible for all $q\in[1,\infty]$ and $\AAA^{-1}\in\Sigma^\wegt_{\alpha^{-1}}$.
Moreover, if $\AAA\in\E_{\alpha}$, then we also have $\AAA^{-1}\in\E_{\alpha^{-1}}$. 
\et
\bpf
First, we observe that $\AAA^{-1} = (\AAA^\star\AAA)^{-1}\AAA^\star$.
Second, since Lemma \ref{bp} implies $\AAA^\star\AAA\in\Sigma^\wegt_1$ (or $\E_1$),
\cite[Theorem 2]{Bas97'} guarantees that $ (\AAA^\star\AAA)^{-1}\in\Sigma^\wegt_1$ (or $\E_1$). Finally, applying Lemma 
\ref{bp} once again we get the desired results.
\epf

\begin{rem}\label{wien}
The above result may seem remarkable but Lemma \ref {PL} shows that in most interesting cases  
it is vacuous unless $|\alpha| = 1$.
The case $\alpha = 1$, however, is standard and to prove the result when $\alpha = -1$ it is enough to recall that if $\AAA\in\Sigma_{-1}$ then 
$\AAA^2\in\Sigma_1$ by Lemma \ref{bp} or, if $\AAA^2$ is not well-defined, employ the reindexing trick used in Remark \ref{R123}. 
The following is a less trivial extension of Wiener's Lemma.
\end{rem}

\bt
\label{corHilbert}
Let $X_n = \mathcal H_X$ and $Y_n = \mathcal H_Y$ be the same Hilbert (or Euclidean) spaces for all $n\in\GG$  and $\AAA\in\Sigma^\wegt_\alpha$ where $\wegt(j)=(1+|j|)^s$, $ s>(d+1)^2$. Let also $p\in[1,\infty]$.
\begin {enumerate}
\item [(i)] If $\AAA$ is $p$-bb, then $\AAA$ is left invertible for all $q\in[1,\infty]$ and a left inverse is given by  
$\AAA^{\sharp}=(\AAA^\star\AAA)^{-1}\AAA^\star\in\Sigma^\wegt_{\alpha^{-1}}$.
\item [(ii)] If $\AAA^\star$ is $p$-bb, then $\AAA$ is right invertible for all $q\in[1,\infty]$ and a right inverse is given by  $\AAA^{\flat} =\AAA^\star(\AAA\AAA^\star)^{-1}\in\Sigma^\wegt_{\alpha^{-1}}$.
\end {enumerate}  
\et
\bpf Since (i) and (ii) are equivalent, we prove only (i).
Theorem \ref{mainthm} implies that $\norm{\AAA x}_2\ge \kappa_2 \|x\|_2$ for some $\kappa_2 > 0$ and all $x\in\ell^2$. 
Under the specified conditions the Banach spaces $\ell^2(\GG,(X_n))$ and $\ell^2(\GG,(Y_n))$ are, however,
Hilbert spaces and $\AAA^\star$ defines the Hilbert adjoint of $\AAA$. Since $\la\AAA^\star\AAA x, x \ra = \la\AAA x, \AAA x \ra \geq \kappa_2
 \la x, x \ra$, we have that the operator $\AAA^\star\AAA$ is invertible in $\ell^2$. It remains to argue as in Theorem \ref{Inv} and apply Lemma \ref{bp} and \cite[Theorem 2 and Corollary 3]{Bas97'}. 
 \epf
\bc\label{corcompl}
If $\AAA$ is as in Theorem \ref{corHilbert}(i) then $Im\,\AAA$ is a complementable subspace of $\ell^q$, $q\in[1,\infty]$.
\ec

\section{Applications}
\label {A}

This section is mainly devoted to two applications. The first application concerns Banach frames and the second one concerns sampling theory.

\subsection{Banach Frames.}\
\label {BF}
The notion of a frame in a separable Hilbert space has already become classical. The pioneering work \cite{DS} explicitly introducing it
 was published in 1952. Its analogues in Banach spaces, however, are non-trivial (see \cite{AST, BCHL, C 03, Gro} and references therein). In this subsection we show that in case of
certain localized frames  the simplest possible extension of the definition remains meaningful.

\begin{defn} 
Let $\HH$ be a separable Hilbert space. A sequence $\varphi_n\in\HH$, $n\in\GG$, is a \textit{frame} for $\HH$ if
for some $0<a\le b<\infty$
\begeq\label{eqnorm}
a\norm{f}^2 \le \sum_{n\in\GG} |\la f,\varphi_n \ra|^2 \le b\norm{f}^2
\eq
for all $f\in\HH$.
\end{defn}

The operator $T: \HH\to \lt$, $Tf = \{\la f,\varphi_n \ra\}_{n\in\GG}$, $f\in\HH$, is called an \textit{analysis operator}.
It is an easy exercise to show that a sequence $\varphi_n\in\HH$ is a frame for $\HH$ if and only if its analysis operator has a left inverse.
The adjoint of the analysis operator, $T^*:\lt\to\HH$, is given by $T^*c = \sum\limits_{n\in\GG} c_n\varphi_n$, $c = (c_n)\in\lt$. The \textit{frame operator} is
$T^*T: \HH\to\HH$, $T^*Tf = \sum\limits_{n\in\GG} \la f,\varphi_n \ra\varphi_n$, $f\in\HH$. Again,
a sequence $\varphi_n\in\HH$ is a frame for $\HH$ if and only if its frame operator is invertible. The canonical dual frame $\tilde\varphi_n\in\HH$ 
is then $\tilde\varphi_n = (T^*T)^{-1}\varphi_n$ and the (canonical) \textit{synthesis operator} is
$T^{\sharp}: \lt\to H$, $T^{\sharp} = (T^*T)^{-1}T^*$, so that
\[f = T^{\sharp} Tf = \sum_{n\in\GG} \la f,\varphi_n \ra\tilde\varphi_n = \sum_{n\in\GG} \la f,\tilde\varphi_n \ra\varphi_n\]
for all $f\in\HH$.

In general Banach spaces one cannot use just the equivalence of norms similar to \eqref{eqnorm}. The above construction breaks down
because, in this case, the analysis operator ends up being bounded below and not necessarily left invertible. As a result a ``frame
decomposition'' remains possible but ``frame reconstruction'' no longer makes sense. Theorem \ref{corHilbert}(i) indicates, however,
that often this obstruction does not exist. The idea of this section is to make the previous statement precise. To simplify the
exposition we remain in the realm of Banach spaces $\ell^p(\GG,\HH)$ and use other chains of spaces such as as the one in \cite{Gro}
only implicitly.

\begin{defn}\label{pframe}
 A sequence $\varphi^n = (\varphi^n_m)_{m\in\GG} \in \ell^1(\GG,\HH)$, $n\in\GG$, is a  $p$-{\it frame} (for $\ell^p(\GG,\HH)$) for some $p\in[1,\infty)$ if
\begeq\label{peqnorm}
a\norm{f}^p \le \sum_{n\in\GG} \left\vert\sum_{m\in\GG}\la f_m,\varphi_m^n \ra\right\vert^p \le b\norm{f}^p
\eq
for some $0<a\le b<\infty$ and all $f = (f_m)_{m\in\GG}\in\ell^p(\GG,\HH)$. 
If
\begeq\label{peqnorm1}
a\norm{f} \le \sup_{n\in\GG} \left\vert\sum_{m\in\GG}\la f_m,\varphi_m^n \ra\right\vert \le b\norm{f}
\eq
for some $0<a\le b<\infty$ and all $f = (f_m)_{m\in\GG}\in\ell^{\infty}(\GG,\HH)$, then the sequence
$\varphi^n$ is called an $\infty$-frame. It is called a $0$-frame if \eqref{peqnorm1} holds for
all $f \in \ccc_0(\GG,\HH)$. 

\end{defn}

The operator $T_\varphi = T: \lp(\GG,\HH)\to \lp(\GG) = \lp(\GG,\CC)$, given by
\[Tf = \la f,\varphi_n \ra := \{\sum\limits_{m\in\GG}\la f_m,\varphi_m^n \ra\}_{n\in\GG},\quad f\in\lp(\GG,\HH),\]
 is called a \textit{p-analysis operator}, $p\in[1,\infty]$. The $0$-analysis operator is defined the same way for
$f \in \ccc_0(\GG,\HH)$. 

\begin{defn}
A $p$-frame $\varphi^n$ with the $p$-analysis operator $T$, $p\in\{0\}\cup[1,\infty]$, is $(s,\alpha)$-\textit{localized}  for some $s>1$ and $\alpha\ne 0$, if there exists an isomorphism
$J: \linf(\GG,\HH)\to \linf(\GG,\HH)$ which leaves invariant $\ccc_0$ and all $\lqq(\GG,\HH)$, $q\in[1,\infty)$, and such that
\[TJ_{\vert \lp}\in\Sigma_\alpha^\wegt,\]
 where $\wegt(n) = (1 +|n|)^s$, $n\in\GG$, see Remark \ref {OpToMat}.
\end{defn}
\begin {rem}
If $\HH$ is finite dimensional, then the above definition is vacuous for $|\alpha|>1$, due to Lemma \ref {PL}. 
\end {rem}
As a direct corollary of Theorem \ref{corHilbert} and the above definition we obtain the following result.

\bt\label{pfr}
Let $\varphi^n$, $n\in\GG$, be an $(s,\alpha)$-localized $p$-frame  for some $p\in\{0\}\cup[1,\infty]$ with $s> (d+1)^2$. Then 
\begin{enumerate}
\item[(i)] The $q$-analysis operator $T$ is well defined and left invertible for all $q\in\{0\}\cup[1,\infty]$, and the
$q$-synthesis operator $T^{\sharp} = (T^*T)^{-1}T^*$ is also well defined for all $q\in\{0\}\cup[1,\infty]$.

\item[(ii)] The sequence $\varphi^n$, $n\in\GG$, and its dual sequence $\tilde\varphi^n = (T^*T)^{-1}\varphi^n$, $n\in\GG$,
are both $(s,\alpha)$-localized $q$-frames for all $q\in\{0\}\cup[1,\infty]$. 

\item[(iii)] In $\ccc_0$ and $\ell^q$, $q\in[1,\infty)$, we have the reconstruction formula
\[f = T^{\sharp}Tf = \sum_{n\in\GG} \la f,\varphi_n \ra\tilde\varphi_n = \sum_{n\in\GG} \la f,\tilde\varphi_n \ra\varphi_n.\]
For $f\in\linf$ the reconstruction formula remains valid provided the convergence is understood in the weak$^*$-topology.
\end{enumerate}
\et
Theorem \ref {pfr}(iii) shows that an $(s,\alpha)$-localized $p$-frame is a Banach frame for $\ccc_0$ and all $\ell^q$, $q\in [1,\infty]$, in the sense of the following definition.
\begin {defn} \cite [Def. 13.6.1]{grbook} A countable sequence $\{x_n\}_{x_n\in J} \subset X'$ in the dual of a Banach space $X$ is a Banach frame  for $X$ if there exist an associated sequence space $X_d(J)$, a constant $C\ge 1$, and a bounded operator $R: X_d\to X$ such that for all $f \in X$
$$ \frac {1} {C}\|f\|_X\le \|\langle f, x_n\rangle\|_{X_d}\le C \|f\|_X,$$
$$R(\langle f, x_n\rangle_{j\in J})=f.$$
\end {defn}
\begin{exmp}
Following \cite{grbook}, let $g\in\mathcal S\subset C^{\infty} (\RR^{d})$ be  a non-zero window function in the Schwartz class  $\mathcal S$, and $V_g$ be the 
short time Fourier transform 
\[(V_g f)(x,\omega) = \int_{\RR^{2d}} f(t)
\overline{g(t-x)}e^{-2\pi it\cdot\omega} dt,\quad x,\omega\in\RR^d.\]
 Let $ M^p$, $1\le p\le \infty$, be the modulation spaces of tempered distributions with the norms 
$$
\|f\|_{M^p}=\left ( \int_{\RR^d} \left(  \int_{\RR^d} |V_gf(x,\omega)|^pdx \right) d\omega \right)^{1/p}, \quad 1\le p<\infty,
$$

$$
\|f\|_{M^\infty}=\| V_gf\|_\infty.$$
It is known that these modulation spaces do not depend on the choice of $g \in \mathcal S$ and are isomorphic to $\ell^p(\ZZ^{2d})$, with isomorphisms provided by the Wilson bases. 

Let $g\in M^1$ be a window such that the Gabor system
$$
\mathcal G(g,a,b)=\{g_{k,n}(x)= e^{-2\pi i(x-ak)\cdot bn }g(x-ak), \; k,n \in \ZZ^d, x \in \RR^d \}
$$
is a {\it tight Banach frame } for all $M^p$, $1\le p \le \infty$. By this we mean that
the $p$-analysis operator $T_{\mathcal G}: M^p\to \lp$,
$T_{\mathcal G} f = \{\la f, g_{k,n} \ra\}$,  is
left invertible and the frame operator $T^*_{\mathcal G}T_{\mathcal G}$ is a scalar multiple of the identity operator for all $p\in[1,\infty]$. Assume that a sequence $\Phi=\{ \phi_{i,j}\}_{i,j\in\ZZ^d}$ of distributions in $M^\infty$ is such that $\{ \varphi^{(i,j)}=T_{\mathcal G}\phi_{i,j}\}$, $(i,j)\in\ZZ^{2d}$, is an $(s,\alpha)$-localized $p$-frame for some $p \in \{0\}\cup[1,\infty]$, and $s>(d+1)^2$.  Since, by Definition \ref {pframe} $\{ \varphi^{(i,j)}=T_{\mathcal G}\phi_{i,j}\}$ must be in $\ell^1(\ZZ^{2d},\CC)$, then by \cite [Corollary 12.2.8] {grbook}  $\Phi\subset M^1$. Moreover, by Theorem \ref {pfr} we have that $\{ \varphi^{(i,j)}\}$ is an $(s,\alpha)$-localized $q$-frame for all $q \in \{0\}\cup[1,\infty]$, and a Banach frame.  Finally, since $\mathcal G$ is a tight Banach frame for all $M^q$, $q\in [1,\infty]$, we have that 
$$
\langle f,\phi_{i,j}\rangle= Const\langle T^*_{\mathcal G}T_{\mathcal G}f,\phi_{i,j}\rangle=Const\langle T_{\mathcal G}f,T_{\mathcal G}\phi_{i,j}\rangle, \quad \mbox {for all } f \in M^q, 
$$  
and, hence, the frame operator 
$$ f\mapsto T_{\mathcal G}f \mapsto \{\langle T_{\mathcal G}f, T_{\mathcal G}\phi_{i,j}\rangle\} \mapsto\{\langle f, \phi_{i,j}\rangle\}: M^q\to \ell^q(\ZZ^{2d},\CC) $$  is left invertible and, therefore, $\Phi$ is a Banach frame for all $M^q$, $q \in [1, \infty]$. 
\end{exmp}
\begin{rem} A similar example can be produced for more general modulation spaces  $M^{p,q}_\nu$ \cite [Chapters 11, 12] {grbook}. Moreover, the Gabor frame $\mathcal G$ can be replaced by any frame in $M^1$, e.g., a Wilson Basis \cite[Section 12.3]{grbook}.   \end {rem}
\begin{exmp}
Here we would like to highlight the role of the slant $\alpha$ in the previous example. Using the same notation as above, let $\Phi$ be
the frame consisting of two copies of the frame $\mathcal G$. Then
(renumbering $\Phi$ if needed) it is easy to see that the matrix
$(\la \phi_{i,j}, g_{k,n}\ra)_{(i,j),(k,n)\in\ZZ^{2d}}$ is $\frac12$-slanted.
Hence, the slant $\alpha$ serves as a measure of relative redundancy of $\Phi$ with respect to $\mathcal G$ and a measure of absolute redundancy of $\Phi$ if $\mathcal G$ is a basis.
\end{exmp}

\begin{rem}
In the theory of localized frames introduced by K. Gr\"ochenig \cite{Gro}  it is possible to extend a localized (Hilbert) frame to Banach frames for
the associated Banach spaces. The technique we developed in Section 2, allows us to start with a localized $p$-frame and deduce that it is, in fact, a Banach frame for the associated Banach spaces.
Slanted matrices provide us with additional information which makes it
possible to shift emphasis from the frame operator $T^*T$ to the
analysis operator $T$ itself.
\end{rem}

\subsection{Sampling and Reconstruction Problems}\
\label {SRP}
In this subsection we apply the previous results to handle certain problems in sampling theory.
Theorem \ref{mainappl} below was the principal motivation for us to prove Theorem \ref{mainthm}. 

The sampling/reconstruction problem includes devising efficient methods for
representing a signal (function) in terms of a discrete (finite or countable) set of its samples (values) and
reconstructing the original signal from its samples. 
In this paper we assume that the signal is a function $f$ that belongs to a space
\[V^p(\Phi) = \left\{\sum_{k\in \GG} c_k\varphi_k\right\},\]
where $c = (c_k)\in\ell^p(\GG)$ when $p\in[1,\infty]$, $c\in\ccc_0$ when $p =0$,
and $\Phi = \{\varphi_k\}_{k\in\GG}\subset L^p(\RR^d)$ is a countable collection of continuous functions.
To avoid convergence issues in the definition of $V^p(\Phi)$,
we assume that the functions in $\Phi$  satisfy the condition
\begin {equation}
\label {RBC}
m_p\|c\|_{\ell^p}\le \norm{\sum_{k\in\GG} c_k\varphi_k}_{L^p}\le M_p\|c\|_{\ell^p}, \quad \mbox{for all }\, c\in \ell^p,
\end {equation}
for some $m_p,M_p>0$ independent of $c$. This is a  p-Riesz basis condition for $p\in [1,\infty)\cup \{0\}$  \cite {AST}. Furthermore, we assume that the functions in $\Phi$ belong to a Wiener-amalgam space $W^1_\wegt$ defined as follows.
\begin{defn} 
A measurable function $\varphi$  belongs to $W^1_\wegt$ for a certain weight $\wegt$, if it satisfies
\begeq
      \label{cc2}
      \norm {\varphi}_{W^1_\wegt} = \left( \sumzd{k}\wegt(k)\cdot {\ssup \{ \abs {\varphi (x+k)}:\ x
\in \Id \}} \right) < \infty.
\end{equation}
\end{defn}
When a function $\varphi$ in $W^1_\wegt$ is continuous we write $\varphi \in W^1_{0,\wegt}$.
In many applications $V^p(\Phi)$ is a shift invariant space, that is, $\varphi_k(x) = \varphi(x-k)$, $k\in\GG$,
for some $\varphi\in W^1_\wegt$.

Sampling is assumed to be performed by a countable collection of finite complex Borel measures {\Large $\mu$}$ = \{\mu_j\}_{j\in\GG}\subset\mathcal M(\RR^d)$. A {\Large$\mu$}-\textit{sample} is a sequence $f(${\Large $\mu$}$) = \int f d\mu_j$, $j\in\GG$. If $f(${\Large$\mu$}$) \in \ell^p$ and $\|f(${\Large$\mu$}$)\|_{\ell^p}\le C\|f\|_{L^p}$ for all  $f\in V^p(\Phi)$, we say that
{\Large $\mu$} is a $(\Phi,p)$-\textit{sampler}. If a sampler {\Large $\mu$} is a collection of Dirac measures then it is called an $(\Phi,p)$-\textit{ideal} sampler. Otherwise, it is
an $(\Phi,p)$-\textit{average} sampler. 

One of the
main goals of sampling theory is to determine when a sampler {\Large $\mu$} is \textit{stable}, that is when $f$ is uniquely determined
by its {\Large $\mu$}-sample and a small perturbation of the sampler
 results in a small perturbation of $f\in V^p(\Phi)$.
The above condition can be formulated as follows \cite{aldgr:siam}:
\begin{defn} A sampler {\Large $\mu$} 
is \textit{stable on $V^p(\Phi)$ } (in other words,
{\Large $\mu$} is  a stable $(\Phi,p)$-sampler) if the bi-infinite matrix $\AAA^\Phi_\mu$  defined by 
\[(\AAA^\Phi_\mu c)(j) = \sumzd k \int c_k\varphi_k d\mu_j,\ c\in\lp(\GG),\]
defines a bounded  \textit{sampling operator} $\AAA^\Phi_\mu: \lp(\GG)\to\lp(\GG)$ which is bounded below in $\lp$ (or $p$-bb). 
\end{defn}

We assume that the generator $\Phi$ and the sampler {\Large $\mu$} are such that  the operator $\AAA^\Phi_\mu$ is bounded on $\ccc_0$ and all $\lp$, $p\in [1,\infty]$; we say that such sampling system ($\Phi$,{\Large $\mu$}) is \textit{sparse}. This situation happens, for example, when the generator $\Phi$ has sufficient decay at $\infty$ and the sampler is \textit{ separated}. The following theorem is a direct corollary of Theorem \ref{corHilbert} and the above definitions.

\bt\label{mainappl}
Assume that $\wegt(n) = (1 +|n|)^s$, $n\in\GG$, $s > (d+1)^2$, $\Phi$ satisfies \eqref{RBC} for all $q\in\{0\}\cup[1,\infty]$, and {\Large $\mu$} is a $(\Phi,p)$-sampler for every $p \in [1,\infty]$. 
Assume also that the sampling operator $\AAA^\Phi_\mu$ is $p$-bb for some $p\in\{0\}\cup[1,\infty]$ and $\AAA^\Phi_\mu\in\Sigma_\alpha^\wegt$ for some 
$\alpha \neq 0$. 
Then {\Large $\mu$} is a stable sampler on $V^q(\Phi)$ for every $q\in\{0\}\cup[1,\infty]$.
\et

Below we study the case of ideal sampling in shift invariant spaces in greater detail and obtain specific examples of the use of the above theorem.
From now on we assume that $\varphi_k(x) = \varphi(x-k)$, $k\in\GG$, for some $\varphi\in C\cap W^1_\wegt=:W^1_{0,\wegt}$.

\begin{defn} If {\Large $\mu$}$=(\mu_j)$ is a stable ideal sampler on $V^p(\Phi)$ and the measures $\mu_j$ are supported on $\{x_j\}$, $j\in\GG$,
then the set $X = \{x_j,\ j\in\GG\}$ is called a (stable) \textit{set of sampling} on $V^p(\Phi)$.
A set of sampling $X\subset \RR^d$ is {\it separated} if \[\inf_{j\neq k \in\GG} |x_j - x_k| =\delta >0.\]
A set of sampling $X\subset \RR^d$ is {\it homogeneous} if \[\#\{X\cap [n,n+1)\} = M\]
is constant for every $n\in\GG$.
\end{defn}

We are interested in the homogeneous sets of sampling because of the following result.

\bl\label{wegts}
Let  $\varphi\in W^1_{0,\wegt}$, $\Phi = \{\varphi(\cdot -k)\}$, and {\Large $\mu$}$\in \linf(\GG,\mathcal M(\RR^d))$ be an ideal sampler with
a separated homogeneous sampling set $X$. 
Then the sampling operator $\AAA^\Phi_\mu$ belongs to $\Sigma_\alpha^\wegt$ for  $\alpha=M^{-1}$.
\el

\bpf Follows by direct computation. \epf

The following lemma shows that we can restrict our attention to homogeneous sets of sampling
without any loss of generality. The intuition behind this result is that we can count each measurement at a point
in $X$ not once but finitely many times and still obtain unique and stable reconstructions.

\bl\label{homog}
Let $\AAA$ be an infinite matrix that defines a bounded operator on $\ell^p$, $p\in[1,\infty]$, and $\tilde\AAA$ be a (bounded)
operator on $\ell^p$ obtained from $\AAA$ by duplicating each row at most $M$ times. Then $\AAA$ is $p$-bb
if and only if $\tilde\AAA$ is $p$-bb.
\el

\bpf The proof for $p<\infty$ follows from the inequalities
\[\norm{\AAA x}_p^p \leq \|{\tilde\AAA}x \|_p^p \leq (M+1)\norm{\AAA x}_p^p,\ x\in\lp.\]
For $p=\infty$, we have $\norm{\AAA x}_\infty = \|{\tilde\AAA} x\|_\infty$, $x\in\ell^\infty$.
\epf

As a direct corollary of Theorems \ref{corHilbert}, \ref{mainappl},  Lemmas \ref{wegts}, \ref{homog}, and Remark \ref {R22} we obtain the following theorem.

\bt\label{idealsmpl}
Let $\wegt(n) = (1 +|n|)^s$, $n\in\GG$, $s > (d+1)^2$, $\varphi\in W^1_{0,\wegt}$, and
\[a_p\|f\|_{\Lp} \le \|\{f(x_j)\}\|_{\lp} \le b_p\|f\|_{L^p}, \mbox{  for all } f \in V^p(\Phi),\]
for some $p\in [1,\infty]\cup \{0\}$ and a separated set $X = \{x_j,\ j\in\GG\}$. Then $X$ is a stable set of sampling on $V^q(\Phi)$ for all $q\in[1,\infty]\cup\{0\}$.
\et

Now we can prove a Beurling-Landau type theorem \cite{AAK, AG, aldgr:siam, AK} for shift-invariant spaces generated by piecewise differentiable functions.

\bt\label{lowbound}
Let $\Phi$ be a sequence generated by the translates of a piecewise differentiable function $\varphi \in W^1_{0,\wegt}$ such that \[a\|c\|_\infty \leq \norm{\sum_{k\in\ZZ} c_k\varphi_k}_\infty\leq b\|c\|_\infty\quad \mbox{and}\quad
\norm{\sum_{k\in\GG} c_k\varphi^\prime_k}_\infty\leq b^\prime\|c\|_\infty,\quad 
\]
for all $c\in \ccc_0(\GG,(X_n))$. Then every $X= \{x_j\}$ that satisfies $\gamma(X) = \sup(x_{j+1}-x_j) < \frac{2a}{b^\prime}$ is
a set of sampling for $V^p(\Phi)$ for all $p\in\{0\}\cup[1,\infty]$. 
\et

\bpf
We prove the result for everywhere differentiable functions $\varphi$ and omit the obvious generalization.

Let $f\in V^0(\Phi)$ be such that $f^\prime = \sum_{k\in\GG} c_k\varphi^\prime_k$, where the series has  finitely many nonzero terms. The set of such
functions is dense in $V^0(\Phi)$ and if we prove that for all such $f$
\[\norm{\{f(x_j)\}}_\infty = \sup_{j\in \GG}|f(x_j)|\geq \kappa_\infty \|c\|_\infty,\]
the result would follow immediately from Theorem \ref{idealsmpl}.

Let $x^*\in\RR$ be such that $\|f\|_\infty = |f(x^*)|$. 
There exists $j\in J$ such that $|x_j - x^*| \leq \frac12 \gamma(X)$.  
Using the Fundamental Theorem of Calculus, we get 
\[|f(x_j)| = \left\vert\int_{x_j}^{x^*}f^\prime(t)dt - f(x^*)\right\vert\geq
\|f\|_\infty - \left\vert\int_{x_j}^{x^*}\sum_{k\in\GG} c_k\varphi^\prime_k(t)dt\right\vert\]
\[\geq \|f\|_\infty - \left\vert\int_{x_j}^{x^*}\norm{\sum_{k\in\GG} c_k\varphi^\prime_k}_\infty dt\right\vert
\geq (a - \frac12 b^\prime\gamma(X))\|c\|_\infty.\]
Since $\gamma(X) < \frac{2a}{b^\prime}$, we have $\kappa_\infty >a - \frac12 b^\prime\cdot\frac{2a}{b^\prime} = 0$. \epf

\begin{cor}\label{corlowbound}
Let $\Phi$ be a sequence generated by the translates of a piecewise
twice differentiable function $\varphi \in W^1_\wegt$ such that \[a\|c\|_\infty \leq \norm{\sum_{k\in\GG} c_k\varphi_k}_\infty\leq b\|c\|_\infty\quad \mbox{and}\quad
\norm{\sum_{k\in\GG} c_k\varphi^{\prime\prime}_k}_\infty\leq b^{\prime\prime}\|c\|_\infty,\quad 
\]
for all $c\in \ccc_0(\ZZ,(X_n))$. Then every $X= \{x_j\}$ that satisfies $\gamma(X) = \sup(x_{j+1}-x_j) < \sqrt{\frac{8a}{b^{\prime\prime}}}$ is
a set of sampling for $V^p(\Phi)$ for all $p\in\{0\}\cup[1,\infty]$.  
\end{cor}

\bpf
Using the same notation as in the proof of the theorem, we see that
$f^\prime(x^*) = 0$ and, therefore,
\[|f(x_j)| = \left\vert\int_{x_j}^{x^*}f^\prime(t)dt - f(x^*)\right\vert \geq \norm{f}_\infty -
\left\vert\int_{x_j}^{x^*}\int_{t}^{x^*}f^{\prime\prime}(u)dudt \right\vert\geq\]
\[ \norm{f}_\infty - \frac{b^{\prime\prime}}{2} |x^*-x_j|^2 \norm{c}_\infty\geq
(a-\frac18 b^{\prime\prime}\gamma^2(X)) \norm{c}_\infty. \]
At this point the statement easily follows.
\epf

In the next two examples we apply the above theorem and its corollary
to spaces generated by $B$-splines
$\beta_1 =\chi_{[0,1]}*\chi_{[0,1]}$ and $\beta_2 = \chi_{[0,1]}*\chi_{[0,1]}*\chi_{[0,1]}$.

\begin{exmp} Let $\varphi = \beta_1$. This function satisfies the conditions of 
Theorem \ref{lowbound} with $a = 1$ and $b^\prime = 2$. Hence, if 
$\gamma(X) < 1$, we have that $X$ is a set of sampling for $V^0(\varphi)$
with the lower bound $1-\gamma(X)$. Using the estimates in the proof of Theorem \ref{mainthm} one can obtain
explicit lower bounds for any $V^p(\varphi)$, $p\in[1,\infty]$, and a universal bound for all $p\in[1,\infty]$ (see Remark \ref {lb}).
\end{exmp}

\begin{exmp} Let $\varphi = \beta_2$. This function satisfies the conditions of 
Corollary \ref{lowbound} with $a = \frac1{2}$ and $b^{\prime\prime} = 4$. Hence, if 
$\gamma(X) < 1$, we have that $X$ is a set of sampling for $V^0(\varphi)$
with the lower bound $\frac{1}{2}(1-\gamma^2(X))$. Again, using the estimates in the proof of Theorem \ref{mainthm}, one can obtain
explicit lower bounds for any $V^p(\varphi)$, $p\in[1,\infty]$, and a universal bound for all $p\in[1,\infty]$.
\end{exmp}

\subsection {Other Applications.}\ Slanted matrices have also been studied in wavelet  theory and signal processing (see e.g. the book of Bratteli  and Jorgensen \cite {BJ02},  the papers \cite {BP99,BM94,FGR, KDG}, and the references therein). 
In signal processing and communication, a sequence $s$ (a discrete signal) is often split into a finite set of {\it compressed sequences} $\{s_1,\ldots,s_r\}$ from which the original sequence $s$ can be reconstructed or approximated. The compression is often performed with {\it filter banks} \cite {FGR, KDG} which correspond to a set of bi-infinite slanted matrices $\AAA_1,\ldots,\AAA_r$ with corresponding slants $|\alpha_i|\le 1$. Left invertibility and boundedness below is an important issue in these systems and our theory may be useful in such applications. Slanted matrices also occur in $K$-theory of operator algebras and its applications to topology of manifolds \cite{ GY06}. Our results and technique may be applied to these situations as well. Finally, our results may  be useful in the study  of differential equations with unbounded operator coefficients similar to the ones described in \cite{Bas99, BasKri, BasPas}.

\section{Acknowledgments}

First and foremost we would like to thank K. Gr\"ochenig for his comments on \cite{AK} which inspired us to embark on this project.
Secondly, we would like to thank all those people that attended our talks on the above results and shared their valuable opinions, 
to name just a few, R. Balan, C. Heil, P. Jorgensen, G. Pfander, R. Tessera. 
Finally, we thank the cat Rosie for gracefully allowing us to divert our attention from him and type this paper.

\end {document}